\documentclass[a4paper]{article}
\usepackage{amscd,amssymb,amsmath,amsthm}
\usepackage{endnotes} 
\usepackage{mathrsfs}

\makeatletter
\renewcommand{\@seccntformat}[1]
{{\csname the#1\endcsname}.\hspace{0.3em}}
\makeatother

\theoremstyle{plain}
\newtheorem{theorem}{Theorem}
\newtheorem{lemma}[theorem]{Lemma}
\newtheorem{corollary}[theorem]{Corollary}
\newtheorem*{GT*}{Gottlieb's theorem}
\newtheorem*{MT*}{Moser isotopy theorem}
\newtheorem*{CT*}{Compactness theorem}
\newtheorem*{LT*}{Lefschetz-Gottlieb theorem}
\newtheorem*{GSprime*}{Theorem~\ref{GS}$^\prime$}

\theoremstyle{definition}
\newtheorem*{defin*}{Definition}
\newtheorem{example}{Example}
\newtheorem{remark}[example]{Remark}

\DeclareMathAlphabet{\matheur}{U}{eur}{m}{n}
\DeclareMathAlphabet{\matheus}{U}{eus}{m}{n}
\DeclareMathAlphabet{\matheuf}{U}{euf}{m}{n}

\numberwithin{equation}{section}

\newcommand{\norm}[1]{\left\lVert#1\right\rVert}
\newcommand{\pt}{\mathit{pt}}
\newcommand{\diff}{\mathit{Diff}}
\newcommand{\symp}{\mathit{Symp}}

\DeclareMathOperator{\ev}{ev}
\DeclareMathOperator{\Id}{Id}

\DeclareMathOperator{\ind}{ind}
\DeclareMathOperator{\Ker}{Ker}
\DeclareMathOperator{\PSL}{PSL}
\DeclareMathOperator{\PD}{PD}

\begin{document}

\author{Gerasim  Kokarev\thanks{Supported by EPSRC grant GR/S68712/01}
\\ {\small\it School of Mathematics, The University of Edinburgh}
\\ {\small\it King's Buildings, Mayfield Road, Edinburgh EH9 3JZ, UK}
\\ {\small\it Email: {\tt G.Kokarev@ed.ac.uk}}
}

\title{On the topology of the evaluation map and rational curves}
\date{}
\maketitle

\begin{abstract}
We explore a relationship between topological properties of orbits of
$2$-cycles in the symplectomorphism group $\symp (M)$ and the
existence of rational curves in $M$. Under the absence of rational
curves hypothesis, we show that evaluation map vanishies on $\pi_2$ and
obtain a Gottlieb-type vanishing theorem for toroidal cycles in $\symp
(M)$.
\end{abstract}




\section{Statement and discussion of results}
\label{Intro}
\subsection{}
Let $M$ be a closed connected smooth\footnote{by this we always mean
  $C^\infty$-category unless there is an explicit statement to the
  contrary} manifold of dimension $2n$. Suppose $M$ admits a
  symplectic structure $\omega$ and denote by $\mathcal J_\omega$ the
  space of almost complex structures $J$ on $M$ such that $\omega$ is
  $J$-invariant and tames $J$. For an almost complex structure
  $J\in\mathcal J_\omega$ by {\it $J$-curve} in $M$ we mean a
  $J$-holomorphic map $u:\Sigma\to M$, where $\Sigma$ is a closed
  Riemannian surface; $J$-holomorphic spheres are also called {\it
  rational $J$-curves}.  

Recall that a diffeomorphism $\phi:M\to M$ is called {\it
  symplectomorphism} if it preserves the symplectic structure,
  $\phi^*\omega=\omega$. Let $\symp(M)$ be a group of
  symplectomorphisms of $M$ endowed with the compact open topology. By
  $\ev_u$ we denote the evaluation map at a base point $u\in M$ given
  by
$$
\symp(M)\ni\varphi\longmapsto \varphi(u)\in M.
$$
We also use the notation $\ev_u^\natural[\phi]$ for the homotopy class
$[\ev_u\circ\phi]$, the image under the evaluation map of a homotopy
class $[\phi]$ of maps $\Sigma\to\symp(M)$.

Our principle result shows that the existence of a certain $2$-cycle
in $\symp(M)$ whose orbits are non-contractible is related to the
presence of rational $J$-curves.
\begin{theorem}
\label{RaCur}
Let $(M,\omega)$ be a closed symplectic manifold and suppose that
either:
\begin{itemize}
\item[(i)] there exists a homotopy class $[\phi]$ of $2$-spheres in
  $\symp(M)$ such that the evaluated class $\ev_u^\natural[\phi]$ is
  non-trivial;
\item[(ii)] the Euler-Poincar\'e number $\chi(M)$ does not vanish
  and there exists a homotopy class $[\phi]$ of $2$-tori in
  $\symp(M)$ such that the evaluated class $\ev_u^\natural[\phi]$ is
  non-trivial.
\end{itemize}
Then for any almost complex structure $J\in\mathcal J_\omega$ there
exists a $J$-holomorphic sphere in $M$.
\end{theorem}
The proof is based on the analysis of the Cauchy-Riemann equations
perturbed by a term defined by a $2$-parameter family of
symplectomorphisms. We show that the corresponding moduli space
formed by solutions in $\ev_u^\natural[\phi]$ is diffeomorphic to
$M$. In the case~(i) we use the absence of sphere bubbles to conclude
that this moduli space is null-cobordant or the homotopy class
$\ev_u^\natural [\phi]$ is trivial. Further analysis shows that only
the latter occurs. The argument in the case~(ii) is similar, but uses
the Morse-Bott theory. Theorem~\ref{RaCur} is a consequence of a more
general statement in the next subsection.

An example with a torus shows that the conclusion of the theorem no
longer holds if the condition $\chi(M)\ne 0$ in the part~{(ii)} is
dropped. More precisely, the $2$-torus in $\symp(\mathbb T^2)$
given by the action
$$
\mathbb R^2/\mathbb Z^2\ni (s,t):\mathbb T^2\to\mathbb T^2,
\qquad (x,y)\mapsto (x+s,y+t)
$$
evaluates into the fundamental cycle and, clearly, there are no
non-trivial pseudo-holomorphic spheres in $\mathbb T^2$. 

As a direct consequence we arrive at the following vanishing theorem.
\begin{corollary}
\label{SA}
Suppose that a symplectic manifold $(M,\omega)$ does not admit 
non-trivial J-holomorphic spheres for at least one compatible almost
complex structure $J\in\mathcal J_\omega$. Then:
\begin{itemize}
\item[(i)] the evaluation map $\ev_u:\symp(M)\to M$ induces the
  trivial homomorphism on $\pi_2$;
\item[(ii)] if $\chi(M)\ne 0$ , the image $\ev_u^\natural[\phi]$ of
  any homotopy class $[\phi]$ of $2$-tori in $\symp(M)$ is trivial.
\end{itemize}
\end{corollary}
This statement is interesting from the point of view of Gottlieb's
theory~\cite{Gott1,Gott3}. We discuss this in the next subsection in
more detail. Now we end with simple examples of symplectic manifolds
without $J$-spheres. First, we introduce more notation.

Recall that the energy of a map $u$ from a Riemannian surface
$(\Sigma,i_\Sigma)$ to $(M,J)$, where $J\in\mathcal J_\omega$, is
defined as the integral
$$
E(u)=\int_\Sigma\frac{1}{2}\norm{du(z)}^2d\mathit{Vol}_\Sigma(z),
$$
where $\norm{du(z)}$ is the Hilbert-Schmidt norm of the operator
$du(z):T_z\Sigma\to T_{u(z)}M$ with respect to the metric
$\omega(\cdot,J\cdot)$ on $M$ and any metric on $\Sigma$ in the
conformal class determined by $i_\Sigma$. Any $J$-holomorphic curve
$u:\Sigma\to M$ minimises the energy in a given homology class and, in
particular, enjoys the following identity
\begin{equation}
\label{EnergyId}
E(u)=\left\langle u^*[\omega],\Sigma\right\rangle,
\end{equation}
where the right-hand side stands for the evaluation of $u^*[\omega]$ on
the fundamental cycle. As is known~\cite[Chapter~4]{McDS}, and can be
easily proved, the quantity
$$
S_\omega(J)=\inf\left\{E(u): u\text{ is a non-constant }J\text{-sphere
  in }M\right\}
$$
is positive. Here the infimum over the empty set is supposed to be
equal to infinity. The latter, for example, occurs in the examples
below.
\begin{example}
A symplectic manifold $(M,\omega)$ is called {\it symplectically
  aspherical} if  $\left.\omega\right|_{\pi_2}=0$. As follows from the
  energy identity for $J$-curves (relation~\eqref{EnergyId} above)
  such manifolds do not have non-trivial $J$-spheres for any
  $J\in\mathcal J_\omega$. The existence of symplectically aspherical
  manifolds with non-trivial $\pi_2$ was an open question until the
  examples due to Kollar and Gompf~\cite{Gompf} appeared. There are
  also examples of the latter with arbitrarily large Euler-Poincar\'e
  numbers; it was observed in~\cite{FMS} that one can construct these
  as the symplectic submanifolds described by Auroux~\cite{Aur}. 
\end{example}
\begin{example}
\label{4dim}
Let $(M,\omega)$ be a $4$-dimensional symplectic manifold such that
its first Chern class is a non-positive multiple of $\omega$,
\begin{equation}
\label{SphericalEC}
\left[c_1]\right|_{\pi_2}=k\left.[\omega]\right|_{\pi_2},
\qquad\text{where }k\leqslant 0,\quad k\in\mathbb R. 
\end{equation}
Then for a generic almost complex structure $J\in\mathcal J_\omega$ the
manifold $M$ contains no $J$-spheres. Indeed, first recall that a
$J$-curve is called simple if it is not a (branched) cover of degree
greater than one of another $J$-curve. Clearly, if there exists a
$J$-sphere, then there exists a simple $J$-sphere. Further, due to the
standard Fredholm theory~\cite[Chapter~3]{McDS}, for a generic $J$ the
dimension of unparametrised simple $J$-spheres representing $A\in
H_2(M,\mathbb Z)$ is equal to $2c_1(A)-2$. In particular, if such a
sphere exists, then $c_1(A)\geqslant 1$. On the other hand, due
to relation~\eqref{SphericalEC}, we have $c_1(A)\leqslant 0$ for any 
homology class $A$ that can be represented by a rational $J$-curve.
The latter is a consequence of the energy identity for $J$-curves,
relation~\eqref{EnergyId}.
\end{example}

\subsection{}
\label{energySect}
In this subsection we state a more precise condition on a $2$-cycle
which guarantees that its image under the evaluation map is
homotopically trivial. For this we define a certain energy-type
characteristic of its action on $M$. 

First, we suppose that tori and spheres in $\symp(M)$ under
consideration are represented by maps $\phi$ such that
$\ev_u\circ\phi$ is smooth for any $u\in M$. Due to the following
lemma this does not affect the topological conclusions. 
\begin{lemma}
\label{moser}
Let $\Sigma$ and $(M,\omega)$ be an arbitrary closed manifold and a
closed symplectic manifold respectively. Then continuous maps
$\phi:\Sigma\to\symp(M)$ such that $\ev_u\circ\phi$ is smooth for any
$u\in M$ form a dense subset in the space of all continuous maps
$\Sigma\to\symp(M)$ with respect to the compact open topology.
\end{lemma}
The proof follows essentially from the Moser isotopy argument and is
explained at the end of Section~\ref{proofs}.

Now we define the {\it evaluation energy} of a map $\phi:\Sigma \to
\symp(M)$ as
$$
E^{\ev}(\phi,\omega,J)=\int\limits_\Sigma\frac{1}{2}\max_{u\in
  M}\norm{d(\ev_u\circ\phi)}^2d\mathit{Vol}_\Sigma,
$$
where the Hilbert-Schmidt norm of $d(\ev_u\circ\phi)$ is taken with
respect to the metric $\omega(\cdot,J\cdot)$ on $M$ and some (and,
hence, any) metric in the conformal class of $i_\Sigma$. Further, a map
$\phi:\Sigma \to\symp(M)$ and a given metric on $\Sigma$ define a
function $\Lambda$ on the product $\Sigma\times M$ by the relation
$(\ev_u\circ\phi)^*\omega=\Lambda_u \mathit{dVol}_\Sigma$. Using this
function we construct the second functional 
$$
\Delta(\phi,\omega)=\int_\Sigma\left(\max_u\Lambda_u-\min_u\Lambda_u
\right)d\mathit{Vol}_\Sigma\geqslant 0.
$$
Finally, the {\it corrected evaluation energy} is defined as the sum
$$
\mathrsfs{E}(\phi,\omega,J)=E^{\ev}(\phi,\omega,J)+
\Delta(\phi,\omega).
$$
In general, this quantity depends on the conformal class of metrics or,
equivalently, the complex structure $i_\Sigma$ on the Riemannian
surface. 

Let $\mathcal M_g$ be the Riemannian moduli space of all complex
structures on a Riemannian surface $\Sigma$ of genus $g$ up to the
pull-back by an orientation preserving diffeomorphism. Recall that for
a sphere and a torus the space $\mathcal M_g$ is identified with a
single point and the fundamental domain for the action of
$\PSL(2,\mathbb Z)$ on the upper half-plane, respectively. Denote by
$\mathrsfs{E}_{\Pi}([\phi],\omega,J)$ the infimum of the corrected
evaluation energy $\mathrsfs{E}$ over pairs $(\phi,i_{\Sigma})$, where
$\phi$ represents a given homotopy class $[\phi]$ and $i_{\Sigma}$
ranges over $\mathcal M_g$. 

We are ready to state a quantitative version of Theorem~\ref{RaCur}.
\begin{theorem}
\label{GS}
Let $(M,\omega)$ be a symplectic manifold and $\Sigma$ be a Riemannian
surface. Suppose that a homotopy class $[\phi]$ of mappings
$\Sigma\to\symp(M)$ is such that
\begin{equation}
\label{EnergyC}
\sup_{J\in\mathcal J_\omega}\left(S_\omega(J)-\mathrsfs{E}_\Pi
([\phi],\omega,J)\right)>0.
\end{equation}
Then:
\begin{itemize}
\item[(i)] if $\Sigma$ is a sphere, the homotopy class
  $\ev_u^\natural[\phi]$ is trivial;
\item[(ii)] if $\Sigma$ is a torus, the homotopy class
  $\ev_u^\natural[\phi]$ is trivial or $\chi(M)=0$.
\end{itemize}
\end{theorem}
The statement of the theorem can be also regarded as an estimate for
the energy $\mathrsfs{E}(\phi,\omega,J)$ from below. That is the
``energy'' required for a sphere or a torus in $\symp(M)$ to evaluate
into a homotopically non-trivial one is at least $S_\omega(J)$.

The hypothesis on the evaluation energy in Theorem~\ref{GS} can be
relaxed, if we are concerned only with the action of the evaluation
map on homology classes. Let $A$ be a class from $H_2(\symp,\mathbb
Z)$. Denote by $\mathrsfs{E}_H(A,\omega,J)$ the infimum of the
evaluation energy over pairs $(\phi,i_\Sigma)$, where $\phi$ is a map
of a Riemannian surface $\Sigma$ of a fixed genus $g$ into $\symp(M)$
such that $\phi_*[\Sigma]=A$ and $i_\Sigma\in\mathcal M_g$. We have
the following version of Theorem~\ref{GS}. 
\begin{GSprime*}
Let $(M,\omega)$ be a symplectic manifold and $A$ be a homology class
in $H_2(\symp,\mathbb Z)$ that can be represented by an image of a
given Riemannian surface $\Sigma$. Suppose that
\begin{equation}
\label{EnergyCH}
\sup_{J\in\mathcal J_\omega}\left(S_\omega(J)-\mathrsfs{E}_H
(A,\omega,J)\right)>0.
\end{equation}
Then:
\begin{itemize}
\item[(i)] if $\Sigma$ is a sphere, the homology class $(\ev_u)_* A$
  is trivial;
\item[(ii)] if $\Sigma$ is a torus, the homology class $(\ev_u)_* A$
  is trivial or $\chi(M)=0$.
\end{itemize}
\end{GSprime*}
\begin{remark}
It is a simple exercise to show that the infimums of the corrected
evaluation energy $\mathrsfs{E}$ on the homology classes $A$ and $-A$
coincide. Thus, condition~\eqref{EnergyCH} is natural with respect to
the fact that the map $\ev_u$ vanishes or not on these classes
simultaneously. 
\end{remark}
Finally, we mention that our results can be viewed as symplectic
versions of Gottlieb's vanishing theorems. To illustrate the
relationship more clearly we recall the following assertion, which is
due to~\cite[Theorem~8.9]{Gott2}.
\begin{GT*}
Let $N$ be a closed oriented manifold and $\diff(N)$ be its group of
diffeomorphisms. Then the homomorphisms $\chi(N)\ev^*_u$ and $c\cdot
\sigma(N)\ev_u^*$ of the cohomology groups $H^k(N,R)\to H^k(\diff,R)$ 
vanish for any $k>0$ and any unitary ring $R$; here $c$ is an
appropriate non-zero integer which depends only on the dimension of
$N$ and $\chi(N)$ and $\sigma(N)$ stand for the Euler-Poincar\'e
number and the signature respectively.
\end{GT*}

As the example below shows the Euler-Poincar\'e number $\chi(N)$ in
the theorem is essential and, in general, the homomorphism induced by
$\ev_u$ is not expected to be trivial on cohomology or homology.
\begin{example}
\label{SphereEx}
Let $S^2$ be a unit sphere in $\mathbb R^3$ and $SO(3)$ be its group
of orientation preserving isometries. The evaluation map $\ev_u:SO(3)
\to S^2$ defines a bundle with fibre $SO(2)$. Note that this map
induces the trivial homomorphisms on the reduced homology. However,
the homomorphism on the cohomology 
$$
\ev_u^*:H^2(S^2,\mathbb Z_2)\to 
H^2(SO(3),\mathbb Z_2)
$$ is not trivial. Indeed, the fundamental class $[\omega]=\PD [\pt]$ 
maps to
$$
\ev_u^*[\omega]=\ev_u^*\PD [\pt]=\PD [\ev_u^{-1}(\pt)].
$$
Since the fiber $\ev_u^{-1}(\pt)$ is not homologous to zero and the
Poincar\'e Duality $\PD$ is an isomorphism, we conclude that $\ev_u^*
[\omega]\ne 0$. This illustrates Gottlieb's theorem -- the presence of
the Euler-Poincar\'e number is essential. In particular, we see that
the evaluation map is not contractible on the $2$-skeleton of
$SO(3)$. In fact, there are $2$-tori in $SO(3)$ which evaluate into
homotopically non-trivial ones and, hence, condition~\eqref{EnergyC}
in Theorem~\ref{GS} is necessary. As such a torus one can take, for
example, a subset in $SO(3)$ generated by rotations around two
different axes in $\mathbb R^3$; since $SO(3)$ is not commutative one
needs to specify which rotation applies first.
\end{example}

\paragraph{Acknowledgements.} I am much obliged to Elmer Rees for a
number of discussions on the subject, which encouraged me to write
this note.

\section{Preliminaries}
\subsection{Perturbed Cauchy-Riemann equations}
Let $\Sigma$ be an oriented closed Riemannian surface and $(M,\omega)$
be a closed symplectic manifold of dimension $2n$ endowed with an
almost complex structure $J\in\mathcal J_\omega$. For mappings
$u:\Sigma\to M$ we consider the non-linear Cauchy-Riemann operator
$$
\bar\partial u=\frac{1}{2}(du+J\circ du\circ i_\Sigma),
$$
the $J$-complex anti-linear part of the differential $du$. Denote
by $\Omega^{0,1}$ the vector bundle with base $\Sigma\times M$ whose
fibre over $(z,u)$ is formed by $J$-anti-linear operators $T_z\Sigma\to
T_uM$. In this notation the differential operator $\bar\partial$ sends
$$
\mathit{Maps}(\Sigma,M)\ni u\longmapsto\bar\partial u\in
\mathit{Sections}(\tilde u^*\Omega^{0,1}),
$$
where $\tilde u:\Sigma\to\Sigma\times M$ is the graph of $u$, given by
$z\mapsto (z,u(z))$. More generally, let $\matheur f$ be a section of
the bundle $\Omega^{0,1}$. Consider the {\it perturbed
  Cauchy-Riemann equations}
\begin{equation}
\label{CReq}
\bar\partial u(z)=\matheur f(z,u(z)),\qquad z\in\Sigma;
\end{equation}
its solutions are called {\it perturbed $J$-curves}. Below we suppose
that the right-hand side $\matheur f$ is $W^{p,\ell+1}$-smooth in the
Sobolev sense, where $p>2(n+1)$ and $\ell> 3$, and a solution $u$ is
$W^{2,2}$-smooth. Due to elliptic regularity theory, these
suppositions imply that solutions of equation~\eqref{CReq}  are, in
fact, $C^{\ell+1}$-smooth. 

For a given homotopy class $[v]$ of mappings $\Sigma\to M$ denote by
$\mathfrak M([v],J)$ the {\it universal moduli space} formed by  pairs
$(u,\matheur f)$ of such maps $u\in [v]$ and sections $\matheur f$ of
$\Omega^{0,1}$ which satisfy equation~\eqref{CReq}. We consider
$\mathfrak M([v],J)$ as a subspace in the product
$C^{\ell+1}(\Sigma,M)\times\{C^\ell$-smooth sections $\matheur f\}$ and
endow it with the induced topology. The symbol $\pi$ denotes the
natural projection
$$
\mathfrak M([v],J)\ni (u,\matheur f)\mapsto\matheur
f\in\{W^{p,\ell+1}\text{-smooth sections}\}.
$$
Thus, each fiber $\pi^{-1}(\matheur f)$ is simply the moduli space of
solutions (homotopic to $v$) of equation~\eqref{CReq} with a given
section $\matheur f$.

It is a simple exercise to show that a solution of
equation~\eqref{CReq} satisfies the following energy estimate:
$$
E(u)\leqslant\int_\Sigma\max_u\norm{\matheur
  f(\cdot,u)}^2\mathit{dVol}_{\Sigma}+\langle u^*[\omega],\Sigma\rangle.
$$
Using this and the standard rescaling technique we arrive at
the following statement.
\begin{CT*}
Let $(M,\omega)$ be a closed symplectic manifold endowed with an
almost complex structure $J\in\mathcal J_\omega$. Denote by $\mathrsfs
{C}$ the set formed by homotopy classes $[v]$ of mappings $\Sigma\to M$
such that 
$$
V_{\mathrsfs C}=\sup\left\{\langle v^*[\omega],\Sigma\rangle:
[v]\in\mathrsfs{C}\right\}<S_\omega(J).
$$
Then the natural projection
$$
\pi:\bigcup_{[v]\in\mathrsfs{C}}\mathfrak M([v],J)\to
\{W^{p,\ell+1}\text{-smooth sections }\matheur f\},
\qquad (u,\matheur f)\mapsto\matheur f,
$$
restricted on the domain $\pi^{-1}(\matheus U_\ell)$ is proper, where
\begin{equation}
\label{U}
\matheus U_\ell=\left\{W^{p,\ell+1}\text{-smooth }\matheur f:\int_\Sigma
\max_u\norm{\matheur f(\cdot,u)}^2<S_\omega(J)-V_{\mathrsfs{C}}\right\}.
\end{equation}
In particular, the space of solutions of equation~\eqref{CReq} within
the homotopy classes $[v]$ such that $\langle v^*[\omega],
\Sigma\rangle\leqslant 0$ is always compact provided the section
$\matheur f$ satisfies $\int\max_u\norm{\matheur f(\cdot,u)}^2
<S_\omega(J)$.
\end{CT*}
In applications below the set $\mathrsfs{C}$ is a single homotopy
class or the set of homotopy classes representing a given homology
class of mappings. In both cases the constant $V_{\mathrsfs{C}}$ is
equal to $\langle v^*[\omega],\Sigma\rangle$.

Now we linearise equation~\eqref{CReq} with respect to a linear
connection $\nabla^\Omega$ on the vector bundle $\Omega^{0,1}$.
By definition the corresponding linearised at a ($C^3$-smooth) map $u$
Cauchy-Riemann operator sends a section $\matheur v$ of the pull-back
bundle $u^*TM$ to a section of $\tilde u^*\Omega^{0,1}$,
$$
\matheur v\longmapsto \left(\bar\partial u\right)_*\!\!\matheur v=
\left.\nabla^\Omega_{\partial/\partial t}\right|_{t=0}
\!\!\left(\bar\partial u_t\right);
$$
here $u_t$ is a family of mappings $\Sigma\to M$ such that
$$
\left.u_t\right|_{t=0}=u\quad\text{and}\quad\left.(\partial/\partial
t)\right|_{t=0}u_t=\matheur v.
$$
Such a connection $\nabla^\Omega$ on the vector bundle $\Omega^{0,1}$
can be, for example, built up from a canonical $J$-linear connection
on $M$ and any Levi-Civita connection (of a metric compatible with the
complex structure) on $\Sigma$. More precisely, let $\nabla$ be a
Levi-Civita connection of the metric $g(\cdot,\cdot)=\omega(\cdot,
J\cdot)$. Then the connection $\widetilde \nabla$ given by
$$
\widetilde\nabla_YX=\nabla_YX-\frac{1}{2}J(\nabla_YJ)X,
$$
where $X$ and $Y$ are vector fields on $M$, is $J$-linear. The
corresponding linearised Cauchy-Riemann operator is given by the
formula
$$
(\bar\partial u)_*\matheur v=(\nabla\matheur v)^{0,1}-\frac{1}{2}
J(u)(\nabla_{\matheur v}J)\partial u
$$
and, in particular, does not depend on a connection on $\Sigma$.
Here $\matheur v$ is a vector field along $u$, the symbol
$(\nabla\matheur v)^{0,1}$ stands for the ($J$-)complex anti-linear
part of the form $\nabla\matheur v$, and $\partial u$ is the
$J$-linear part of $du$. For more details we refer
to~\cite[Chapter~3]{McDS}. 

Analogously, the linearisation of equation~\eqref{CReq} at a map $u$
defines the differential operator $(\bar\partial u)_*-\matheur
f_*(\cdot,u)$. This operator differs from the linearised
Cauchy-Riemann operator by zero-order terms depending on derivatives
of $\matheur f$. Moreover, the corresponding operator linearised at a
solution of equation~\eqref{CReq} does not depend on the choice of a
connection $\nabla^\Omega$ used and can be defined as 
$$
\matheur v\longmapsto\left.\frac{\partial}{\partial t}\right|_{t=0}
\left[\bar\partial u_t-\matheur f(\cdot,u_t(\cdot))\right],
$$
where $u_t$ is a family of mappings as above. Recall that a section
$\matheur f$ in the perturbed Cauchy-Riemann equations is called
{\it regular}, if the cokernel of this differential operator is
trivial for any solution $u$ of equation~\eqref{CReq}. In particular, so
is any section $\matheur f$ for which equation~\eqref{CReq} does not
have solutions, i.e. $\pi^{-1}(\matheur f)=\varnothing$. 

The following statement is folkloric and its analogues are proved by
many authors in different frameworks. Our closest references
are~\cite[Chapter~3]{McDS} and~\cite{GK,KoKu2}.
\begin{theorem}[Folklore]
\label{Folk}
Let $(M,\omega)$ be a symplectic manifold of dimension $2n$ endowed
with an almost complex structure $J\in\mathcal J_\omega$ and $\Sigma$
be a closed oriented Riemannian surface (with a fixed complex
structure). Suppose that a given homotopy class $[v]$ of mappings
$\Sigma\to M$ is such that $\langle v^*[\omega],\Sigma\rangle<
S_\omega(J)$ and let $\matheus U_\ell$ be the domain given
by~\eqref{U} with the integer $\ell$ such that
$$
\ell>n\chi(\Sigma)+2\langle v^*[c_1],\Sigma\rangle+3.
$$
Then for any regular section $\matheur f\in\matheus U_\ell$ the space
of solutions $\pi^{-1}(\matheur f)$ within $[v]$ is either empty or a
closed $C^{\ell-2}$-smooth manifold whose dimension is equal to
$n\chi(\Sigma)+2\langle v^*[c_1],\Sigma\rangle$; besides,
$\pi^{-1}(\matheur f)$ carries a natural orientation. Further, two
regular sections $\matheur f_0$ and $\matheur f_1\in\matheus U_\ell$
can be joined by a path $\matheur f_t\in\matheus U_\ell$ such that
the set $\cup_t\pi^{-1}(\matheur f_t)$ is a smooth oriented manifold
with boundary $\pi^{-1}(\matheur f_0)\cup\pi^{-1}(\matheur f_1)$. The
boundary orientation agrees with the orientation of $\pi^{-1}(\matheur
f_1)$ and is opposite to the orientation of $\pi^{-1}(\matheur f_0)$.
\end{theorem}
We end with a few comments on the proof. First, one shows
that the universal moduli space $\mathfrak M([v],J)$ is a
$C^{\ell-2}$-smooth Banach manifold and the projection $\pi$ is a
$C^{\ell-2}$-smooth Fredholm map. Its index coincides with the
index of the linearised Cauchy-Riemann operator $(\bar\partial u)_*$
and by Riemann-Roch theorem is given by the formula
$$
\ind\pi=n\chi(\Sigma)+2\langle v^*[c_1],\Sigma\rangle.
$$
The regular values of $\pi$ are identified with regular sections
$\matheur f$ and, hence, the preimage $\pi^{-1}(\matheur f)$ is a
$C^{\ell-2}$-smooth manifold whose dimension is equal to $\ind\pi$. 
The proof that two regular fibers are cobordant uses the
transversality argument which requires that the order of smoothness of
$\pi$ is greater than $(\ind\pi+1)$; see also~\cite[Section~3]{Smale}
for a similar argument. This explains the formula for $\ell$ in the
theorem.

For the sequel we point out that the cobordism manifold $N=\cup_t
\pi^{-1}(\matheur f_t)$ is a $C^{\ell-2}$-smooth submanifold of
the universal moduli space $\mathfrak M([v],J)$; the latter is a
submanifold in the product $W^{2,2}(\Sigma,M)\times\{W^{p,\ell+1}
\text{-smooth sections }\matheur f\}$. For given a reference point
$z_*\in\Sigma$ consider the map
\begin{equation}
\label{RefMap}
N\ni (u,\matheur f)\longmapsto u(z_*)\in M.
\end{equation}
The latter factors as the composition of the projection onto
$W^{2,2}(\Sigma,M)$ and the evaluation at the point $z_*$, and 
is clearly $C^{\ell-2}$-smooth.

The case when the dimension of the space of solutions
$\pi^{-1}(\matheur f)$ is equal to zero is of particular interest and
have been studied in~\cite{KoKu2} in a slightly different
framework. We discuss this below in more detail.

\subsection{Elements of Morse-Bott theory}
For the rest of the section we suppose that the genus of a
Riemannian surface $\Sigma$ is equal to one and a given homotopy class
$[v]$ is such that $\langle v^*[c_1],\Sigma\rangle=0$. Then, under the
conditions of Theorem~\ref{Folk}, the space of solutions in $[v]$ of
equation~\eqref{CReq} with a regular $\matheur f\in\matheus U_\ell$ is
finite and its oriented cobordism class defines an integer $\deg\pi$
-- the algebraic number of solutions. Note also that in this case the
linearised Cauchy-Riemann operator sends sections of $u^*TM$ into
themselves (the bundles $u^*TM$ and $\tilde u^*\Omega^{0,1}$ are
naturally isomorphic) and, hence, one can speak about its resolvent
set.

Let $\mathfrak S$ be a space, regarded as a subspace of
$W^{2,2}(\Sigma,M)$, formed by solutions of the equation
\begin{equation}
\label{Grhs}
\bar\partial u(z)=\matheur g(z,u(z)),\qquad z\in\Sigma,
\end{equation}
within a fixed homotopy class. Suppose that $\matheur g$ above is a
smooth section of $\Omega^{0,1}$ and, hence, due to elliptic
regularity, $\mathfrak S$ is formed by smooth mappings. In sequel 
we use the notation $\bar{\mathrsfs{D}} (u)$ for the linearised
operator $(\bar\partial u)_*-\matheur g_*(\cdot,u)$.

By the implicit function theorem any $u\in\mathfrak S$ has a
neighbourhood in the space $\mathfrak S$ which can be identified with
a subset of a ball in the space $\Ker\bar{\mathrsfs{D}}(u)$;
see~\cite[Proposition~4.1]{GK}. In particular, if there exists a
neighbourhood which can be identified with a ball in $\Ker
\bar{\mathrsfs{D}}(u)$, then the space of solutions $\mathfrak S$ is
called {\it non-degenerate at a point} $u$. We call the space
$\mathfrak S$, or its connected component, {\it non-degenerate (in the
  sense of Morse-Bott)} if it is non-degenerate at any point. 
Alternatively, one can say that $\mathfrak S$ is non-degenerate if
each of its connected components $\mathfrak S^\alpha$ is a smooth
submanifold of $W^{2,2}(\Sigma,M)$ whose dimension is equal to the
dimension of $\Ker \bar{\mathrsfs{D}}(u)$, where $u\in\mathfrak
S^\alpha$.
\begin{defin*}
The space of solutions $\mathfrak S$ (or its connected component) is
called {\it strongly non-degenerate} if it is non-degenerate in the
sense of Morse-Bott and for any $u\in\mathfrak S$ the linearised
operator $\bar{\mathrsfs{D}}(u)$ does not have adjoint vectors
corresponding to the zero eigenvalue; i.e. the algebraic multiplicity
of the zero eigenvalue is equal to the dimension of $\Ker
\bar{\mathrsfs{D}}(u)$.
\end{defin*}
\begin{example}
\label{normal}
Suppose that a pull-back bundle $u^*TM$, where $u\in\mathfrak S$, is
endowed with a Riemannian metric. This together with a volume form on
$\Sigma$ gives rise to a natural $L_2$-scalar product on the vector
fields along $u$. Recall that a linear differential operator is called
{\it formally normal} if it commutes with its formally adjoint
operator. Formally normal operators do not have adjoint vectors
corresponding to the zero eigenvalue~\cite[Chapter~5]{Kato}; see 
also~\cite[Section~6.1]{KoKu2}. Thus, if the space of solutions
$\mathfrak S$ is non-degenerate and the operator $\bar{
\mathrsfs{D}}(u)$ is formally normal for any $u\in\mathfrak S$, then
$\mathfrak S$ is strongly non-degenerate.
\end{example}

The following theorem is proved in~\cite{KoKu2};
see~\cite[Theorem~3]{KoKu2} and also the discussion
in~\cite[Section~10]{KoKu2}.
\begin{theorem}
\label{MB}
Let $(M,\omega)$ be a symplectic manifold endowed with an almost
complex structure $J\in\mathcal J_\omega$ and $[v]$ be a homotopy
class of mappings $\Sigma=\mathbb T^2\to M$ such that 
$$
\langle v^*[c_1],\Sigma\rangle=0\quad\text{and}\quad\langle
v^*[\omega],\Sigma\rangle<S_\omega(J).
$$
Suppose that there exists a smooth section $\matheur g\in\matheus 
U_\ell$, $\ell>3$, such that the space $\mathfrak S$ formed by
solutions in $[v]$ of equation~\eqref{Grhs} is strongly non-degenerate
in the sense of Morse-Bott and the evaluation map
\begin{equation}
\label{Sev}
\Sigma\times\mathfrak S\ni (z,u)\stackrel{\ev}{\longmapsto}(z,u(z))
\in\Sigma\times M
\end{equation}
is an embedding. Then the algebraic number $\deg\pi$ of solutions in
$[v]$ of equation~\eqref{CReq} for a regular section $\matheur f\in
\matheus U_\ell$, $\ell>3$, is given by the formula
$$
\deg\pi=\sum_\alpha\pm\chi(\mathfrak S^\alpha),
$$
where $\mathfrak S^\alpha$ is a connected component of $\mathfrak S$
and $\chi(\mathfrak S^\alpha)$ stands for its Euler-Poincar\'e number.
\end{theorem}
\begin{corollary}[Theorem~6 in~\cite{KoKu2}]
\label{MBcor}
Let $(M,\omega)$ be a symplectic manifold endowed with an almost
complex structure $J\in\mathcal J_\omega$. Then the algebraic number
$\deg\pi$ of null-homotopic perturbed $J$-tori for a regular section
$\matheur f$ such that $\int\max_u\norm{\matheur 
f(\cdot,u)}^2<S_\omega(J)$ is equal to the Euler-Poincar\'e number
$\chi(M)$. For a non-trivial homotopy class $[v]$ such that 
$$
\langle v^*[c_1],\Sigma\rangle=0\quad\text{and}\quad\langle
v^*[\omega],\Sigma\rangle \leqslant 0
$$ 
the degree $\deg\pi$ is equal to zero.
\end{corollary}
\begin{proof}
The proof follows directly from Theorem~\ref{MB} by setting $\matheur
g\equiv 0$. Indeed, the space of null-homotopic $J$-tori consists of
all constant mappings only. The corresponding linearised operator 
$\bar{\mathrsfs{D}}(u)$ is the Cauchy-Riemann operator on
vector-functions $\Sigma\to T_uM\simeq\mathbb C^n$. Due to the
Liouville principle $\Ker\bar{\mathrsfs{D}}(u)$ consists of constant
vector-functions only and, hence, the space of null-homotopic
solutions $\mathfrak S\simeq M$ is non-degenerate in the sense of
Morse-Bott. Moreover, the operator $\bar{\mathrsfs{D}}(u)$ is formally
normal and, due to Example~\ref{normal}, we see that $\mathfrak S$ is
strongly non-degenerate. The other hypotheses of the theorem in this
case are obvious. The statement about the non-trivial homotopy class
$[v]$ simply follows from the definition of the degree, since the
suppositions of the theorem imply that $[v]$ does not contain $J$-tori,
i.e. $\pi^{-1}(0)=\varnothing$.
\end{proof}
Note that, since the compactness theorem holds for homology classes of
mappings, Theorem~\ref{MB} also has a version concerned with the
algebraic number of perturbed $J$-tori within homology classes.
In particular, Corollary~\ref{MBcor} implies that for a regular
section $\matheur f$ in equation~\eqref{CReq} such that
$\int\max_u\norm{\matheur f(\cdot,u)}^2<S_\omega(J)$ the algebraic
number of null-homologous perturbed $J$-tori is also equal to
$\chi(M)$. The condition  in Theorem~\ref{MB} that the map given
by~\eqref{Sev} is an embedding can be, in fact, relaxed. 
In~\cite[Appendix~4.B]{GK} it is shown how to deal with the case when
the latter map is only an immersion.

\section{The proofs}
\label{proofs}
Let $\phi:\Sigma\to\diff(M)$ be a fixed map from a Riemannian surface
$\Sigma$ such that the map $(\ev_u\circ\phi)(z)=\phi_z(u)$ is smooth 
with respect to $z\in\Sigma$ for any $u\in M$. Define a section
$\matheur g$ of the bundle $\Omega^{0,1}$ by the following formula:
\begin{equation}
\label{F}
\matheur g(z,u)=\left.\bar\partial (\phi_z(\bar u))
\right|_{\bar u=\phi_z^{-1}(u)}\in\Omega^{0,1}_{(z,u)},
\qquad z\in\Sigma,\quad u\in M. 
\end{equation}
Clearly, for any $u\in M$ the map $\ev_u\circ\phi$ is a solution of
the equation
\begin{equation}
\label{Geq}
\bar\partial u(z)=\matheur g(z,u(z)),\qquad z\in\Sigma.
\end{equation}
Thus, within the homotopy class $\ev_u^\natural[\phi]$ we have the
family of solutions $\{\ev_u\circ\phi\}$ parameterised by $u\in M$.
Our observation is that the Morse-Bott theory applies to
equation~\eqref{Geq}. To implement this we need the following lemmas.
\begin{lemma}
\label{EVid}
For any map $\phi:\Sigma\to\diff(M)$ such that $\ev_u\circ\phi$ is
smooth for any $u\in M$ the following inequality holds:
\begin{equation}
\label{EevFormula}
\int_\Sigma\max_u\norm{\bar\partial(\ev_u\circ\phi)}^2
\mathit{dVol}_\Sigma\leqslant \mathrsfs{E}(\phi,\omega,J)-
\langle (\ev_u\circ\phi)^*[\omega],\Sigma\rangle.
\end{equation}
\end{lemma}
\begin{proof}
Fix a Riemannian metric $g_\Sigma$ within the given conformal class on 
$\Sigma$. Denote by $\Lambda_u$ the function defined by the relation
$(\ev_u\circ\phi)^*\omega=\Lambda_u\mathit{dVol}_\Sigma$. Then 
direct calculations yield the following identity:
$$
\norm{\bar\partial (\ev_u\circ\phi)}^2+\Lambda_u=
\frac{1}{2}\norm{d(\ev_u\circ\phi)}^2.
$$
This implies the inequality
$$
\max_u\norm{\bar\partial (\ev_u\circ\phi)}^2+\min_u\Lambda_u
\leqslant\frac{1}{2}\max_u\norm{d(\ev_u\circ\phi)}^2.
$$
Integrating the latter over $\Sigma$ with respect to the volume form
$\mathit{dVol}_\Sigma$ and using the definition of the functional
$\mathrsfs{E}$ we arrive at the following inequality
$$
\int_\Sigma\max_u\norm{\bar\partial(\ev_u\circ\phi)}^2
\mathit{dVol}_\Sigma\leqslant \mathrsfs{E}(\phi,\omega,J)-
\int_\Sigma\max_u\Lambda_u\mathit{dVol}_\Sigma.
$$
This immediately implies the claim since the last term in the
right-hand side is not greater than $(-\langle 
(\ev_u\circ\phi)^*[\omega],\Sigma\rangle)$. 
\end{proof}

The following lemma is the only place where the hypothesis that $\phi$
takes values in $\symp(M)$ is used.
\begin{lemma}
\label{A}
Suppose that the map $\phi$ takes values in $\symp(M)$. Then
the maps $\ev_u\circ\phi$, where $u\in M$, are the only solutions of
equation~\eqref{Geq} within their homology class, i.e. the space of
mappings $w$ such that $w_*[\Sigma]=(\ev_u\circ\phi)_*[\Sigma]$. In
particular, there are no other solutions within the homotopy class
$\ev_u^\natural[\phi]$.
\end{lemma}
\begin{proof}
Let $w(z)$ be a solution of equation~\eqref{Geq} which is homologous
to the map $(\ev_u\circ\phi)(z)=\phi_z(u)$, where $z\in\Sigma$, and
$u\in M$. Due to the definition of the section $\matheur g$, see
formula~\eqref{F}, this means that
\begin{equation}
\label{reDef}
\bar\partial w(z)=\left.\bar\partial (\phi_z(\bar u))
\right|_{\bar u=\phi_z^{-1}\circ w(z)},\qquad z\in\Sigma.
\end{equation}
Represent the map $w(z)$ as the composition $(\phi_z\circ v)(z)$,
where the map $v:\Sigma\to M$ is defined as $(\phi_z^{-1}\circ
w)(z)$. In particular, the latter map $v$ is null-homologous,
i.e. $v_*[\Sigma]=0$. We need a formula for the value of the
Cauchy-Riemann operator $\bar\partial$ on the composition
$(\phi_z\circ v)(z)$. First, we have
$$
d(\phi_z\circ v)(z)=\left. d_z\phi(\bar u)\right|_{\bar u=v(z)}
+(d\phi_z\circ dv)(z),\qquad z\in\Sigma;
$$
here the left-hand side stands for the differential of the map
$z\mapsto (\phi_z\circ v)(z)$ and by $d_z\phi$ we mean the
differential of the map $\phi_z(u)$ with respect to $z\in\Sigma$. 
Taking the $J$-anti-linear parts of these differentials we arrive at
the following identity
\begin{equation}
\label{keyR}
\bar\partial (\phi_z\circ v)(z)=\left.\bar\partial (\phi_z(\bar u))
\right|_{\bar u=v(z)}+d\phi_z\circ\bar{\Game}v(z),\qquad z\in\Sigma.
\end{equation}
The symbol $\bar{\Game}v$ above denotes the Cauchy-Riemann operator
for the $z$-dependent almost complex structure
$$
\tilde J(z,u)=d\phi^{-1}_z(u)\circ J(\phi_z(u))\circ d\phi_z(u),
\qquad z\in\Sigma,\quad u\in M.
$$
Combining identities~\eqref{reDef} and~\eqref{keyR} and the fact that
$\phi_z$ is a diffeomorphism for any $z\in\Sigma$, we see that the map
$v$ has to satisfy the equation $\bar{\Game}v= 0$. 

For a proof of the lemma we have to show that this map $v$ is
constant; the latter would imply that $w(z)$ is equal to $(\ev_u\circ
\phi)(z)$ for some $u\in M$. The assertion about $v$  follows from an
energy-type identity for solutions of equation $\bar{\Game}v=0$. First, 
it is straightforward to see that for any $z\in\Sigma$ the symplectic
structure $\omega=\phi_z^*\omega$ tames $\tilde J_z$ and is $\tilde
J_z$-invariant. Thus, the bilinear form $\omega(\cdot,\tilde
J_z\cdot)$ defines a scalar product on the tangent space $T_uM$ for
any $z\in\Sigma$ and, hence, a scalar product $\langle\cdot,\cdot
\rangle$ on the space of linear operators $T_z\Sigma\to T_uM$. We
claim that for any solution $v$ of the equation $\bar{\Game}v=0$ the
following energy identity holds
$$
\frac{1}{2}\int_\Sigma\langle dv,dv\rangle\mathit{dVol}_\Sigma
=\langle v^*[\omega],\Sigma\rangle.
$$
This relation immediately implies that a null-homologous solution has
to be constant. Its proof follows the same argument as a proof of the
standart energy identity; see~\cite[Chapter~2]{McDS}.
\end{proof}
\begin{lemma}
\label{B}
The family of solutions $\{\ev_u\circ\phi\}$, where $u\in M$, of
equation~\eqref{Geq} is non-degenerate in the sense of
Morse-Bott. Moreover, if the genus of $\Sigma$ is equal to one, then
this family of solutions is strongly non-degenerate.
\end{lemma}
\begin{proof}
First, we show that the space of solutions $\{\ev_u\circ\phi\}$, where
$u\in M$, is non-degenerate in the sense of Morse-Bott. For this we
have to prove that a vector field $\matheur v$ from the kernel of the
operator $\bar{\mathrsfs{D}}(\ev_u\!\circ\phi)$ has the form $d\phi_z\!
\cdot \matheur v_0$, where $\matheur v_0$ is a vector from
$T_uM$. (Clearly, any vector field of this form belongs to the kernel
of $\bar{\mathrsfs{D}}(\ev_u\circ\phi)$.) Recall that the operator
$\bar{\mathrsfs{D}}(\ev_u\circ\phi)$ can be defined by the following
relation
$$
\bar{\mathrsfs{D}}(\ev_u\circ\phi)\matheur v=\left.
\frac{\partial}{\partial t}\right|_{t=0}\!\!\!\!
\left[\bar\partial w_t-\matheur g(\cdot,w_t(\cdot))\right],
$$
where $w_t$ is a family of mappings $\Sigma\to M$ such that
$$
\left.w_t\right|_{t=0}=\ev_u\circ\phi\quad\text{and}\quad
\left.(\partial/\partial t)\right|_{t=0}w_t=\matheur v.
$$ 
As in the proof of Lemma~\ref{A} we represent $w_t(z)$ as the
composition $(\phi_z\circ v_t)(z)$, where the family of contractible
mappings $v_t(z)$ is defined as $(\phi_z^{-1}\circ w_t)(z)$. In
particular, the map $v_0(z)\equiv u$ is constant. Using
identity~\eqref{keyR} we obtain the following relations
\begin{multline}
\label{keyRlin}
\bar{\mathrsfs{D}}(\ev_u\circ\phi)\matheur v(z)=\left.
\frac{\partial}{\partial t}\right|_{t=0}\!\!
\!\left[d\phi_z\circ\bar{\Game}v_t(z)\right]\\
=d\phi_z\circ\left[\bar{\Game}v_0\right]_*\left.
\frac{\partial}{\partial t}\right|_{t=0}\!\!\!v_t=\left(d\phi_z
\circ\left[\bar{\Game}v_0\right]_*\circ d\phi^{-1}_z\right)
\matheur v(z),\qquad z\in\Sigma.
\end{multline}
This implies that a vector field $\matheur v(z)$ belongs to the
kernel of the operator $\bar{\mathrsfs{D}}(\ev_u\circ\phi)$ if and
only if the composition $d\phi_z^{-1}\cdot\matheur v(z)$, where
$z\in\Sigma$, belongs to the kernel of $[\bar{\Game}v_0]_*$. The
latter operator acts in accordance with the formula
$$
\left[\bar{\Game}v_0\right]_*\!\!\tilde{\matheur v}=\frac{1}{2}
\left(\nabla\tilde{\matheur v}+\tilde J(\cdot,v_0)\circ\nabla
\tilde{\matheur v}\circ i_\Sigma\right)
$$
on sections of the trivial bundle $\mathbb R^{2n}\times\Sigma$ endowed
with the almost complex structure $\tilde J(z,v_0)$ in the fibre over
$z\in\Sigma$. In particular, up to an isomorphism (for example given
by $d\phi_z$) the operator $[\bar{\Game}v_0]_*$ can be regarded as the
usual Cauchy-Riemann operator on vector-functions $\Sigma\to\mathbb
C^{n}$. Due to the Liouville principle, any vector-function from its
kernel has to be constant. Thus, the vector field $d\phi_z^{-1}\cdot
\matheur v(z)$ is constant and we obtain that any $\matheur v$ from
the kernel of $\bar{\mathrsfs{D}}(\ev_u\circ\phi)$ has the form
$d\phi_z\cdot\matheur v_0$ for some vector $\matheur v_0\in
T_uM$. This demonstrates that the space of solutions $\{\ev_u\circ
\phi\}$ is non-degenerate in the sense of Morse-Bott.

Now suppose that $\Sigma$ is a torus. Then the Cauchy-Riemann operator
$[\bar{\Game}v_0]_*$ is formally normal and, in particular, does not
have adjoint vectors corresponding to the zero eigenvalue; see
Example~\ref{normal}. Due to relation~\eqref{keyRlin} so does the
operator $\bar{\mathrsfs{D}}(\ev_u\circ\phi)$. This ends the proof of
the lemma.
\end{proof}
\paragraph{Proof of Theorem~\ref{GS}.}
First, note that the quantities $\mathrsfs{E}(\phi,\omega,J)$ and
$S_\omega(J)$ are invariant under the simultaneous changes
$\omega\mapsto -\omega$ and $J\mapsto -J$, where the map $\phi$ is
arbitrary and an almost complex structure $J$ belongs to $\mathcal
J_\omega$. Besides, the groups of diffeomorphisms preserving the forms
$\omega$ and $(-\omega)$ coincide. Thus, we can suppose that for a
given homotopy class $[\phi]$ the symplectic structure on $M$ is such
that
\begin{equation}
\label{NonPosive}
\langle(\ev_u\circ\phi)^*\omega,\Sigma\rangle\leqslant 0.
\end{equation} 
Under the conditions of the theorem there exist a complex
structure $i_\Sigma$ on $\Sigma$, an almost complex structure
$J\in\mathcal J_\omega$, and a map $\phi:\Sigma\to\symp(M)$,
representing a given homotopy class, such that 
\begin{equation}
\label{ineq}
\mathrsfs{E}(\phi,\omega,J)<S_\omega(J).
\end{equation}
Define a section $\matheur g$ of the bundle $\Omega^{0,1}$ according
to formula~\eqref{F}. Clearly, we have the identity
$$
\max_u\norm{\matheur g(\cdot,u)}^2=\max_u
\norm{\bar\partial(\ev_u\circ\phi)}^2.
$$
Combining this with Lemma~\ref{EVid} and inequality~\eqref{ineq} we 
see that
$$
\int_\Sigma\max_u\norm{\matheur g(\cdot,u)}^2
\mathit{dVol}_\Sigma<S_\omega(J)-\langle(\ev_u\circ\phi)^*\omega,
\Sigma\rangle.
$$
Thus, the section $\matheur g$ belongs to the domain $\matheus U_\ell$
(from the Compactness theorem) given by relation~\eqref{U} with the
constant $V_{\mathrsfs{C}}$ equaled to $\langle(\ev_u\circ\phi)^*
\omega,\Sigma\rangle$. Note that the first Chern class $[c_1](M)$ also
vanishes on the image of $\ev_u\circ\phi$. Indeed, the vector bundle
$(\ev_u\circ\phi)^*TM$ is trivial and, hence, all its characteristic
classes vanish; it is isomorphic to $T_uM\times\Sigma$ under the
morphism which equals $d\phi_z^{-1}(u)$ on the fiber over
$z\in\Sigma$. Hence, due to the Riemann-Roch theorem, the index of the
linearised operator $\bar{\mathrsfs{D}}(\ev_u\circ\phi)$ is equal to
$n\chi(\Sigma)$.

\noindent
{\bf Case~(i).} Suppose $\Sigma$ is a sphere. Then the index of
$\bar{\mathrsfs{D}}(\ev_u\circ\phi)$ is equal to $2n$, the dimension
of $M$. On the other hand, due to Lemmas~\ref{A} and~\ref{B}, 
the space of solutions $\pi^{-1}(\matheur g)$ is formed by the
mappings $\{\ev_u\circ\phi\}$, where $u\in M$, and is non-degenerate
in the sense of Morse-Bott. In particular, we see that $\pi^{-1}(\matheur
g)$ is diffeomorphic to $M$ and the dimension of the kernel of $\bar{
\mathrsfs{D}}(\ev_u\circ\phi)$ is equal to $2n$, the dimension of $M$.
Thus, the index of the operator $\bar{\mathrsfs{D}}(\ev_u\circ\phi)$
is equal to its kernel and, hence, the section $\matheur g$ is
regular. Suppose the homotopy class $\ev_u^\natural[\phi]$ is not
trivial. Then the energy identity~\eqref{EnergyId} and the
hypothesis~\eqref{NonPosive} imply that this homotopy class does not
contain $J$-spheres. Hence, the space of solutions $\pi^{-1}(0)$ is
empty and the zero section of $\Omega^{0,1}$ is also regular as a
right-hand side of equation~\eqref{CReq}. Now Theorem~\ref{Folk}
applies and we see that there is a deformation $\matheur
g_t\in\matheus U_\ell$ of the section $\matheur g$ to zero such that
the preimage $N=\cup_t\pi^{-1}(\matheur g_t)$ is a compact oriented
manifold with boundary $\pi^{-1}(\matheur g)\simeq M$. Choose a
reference point $z_*\in\Sigma$ and consider the map $N\to M$ given by
formula~\eqref{RefMap}. Its restriction to the boundary $\partial N\simeq
M$ acts by the rule
$$
M\ni u\longmapsto \phi_{z_*}(u)\in M
$$
and, in particular, is a diffeomorphism. Thus, we have a
continuous map $N\to\partial N$ whose restriction to the boundary
induces an isomorphism on the top homology. Since the fundamental
class of a closed oriented manifold is non-trivial, the latter is
impossible and, hence, the homotopy class $\ev_u^\natural[\phi]$ has
to be trivial.

\noindent
{\bf Case~(ii).} Suppose $\Sigma$ is a torus. Then the index of
$\bar{\mathrsfs{D}}(\ev_u\circ\phi)$ vanishes and, due to the
discussion in the preceding section, the invariant $\deg\pi$ (the
algebraic number of solutions for a regular section  $\matheur
f\in\matheus U_\ell$) is well-defined. Due to Lemmas~\ref{A} and~\ref{B}
the Morse-Bott theory applies to equation~\eqref{Geq}: the space of
solutions $\pi^{-1}(\matheur g)\simeq M$ is strongly non-degenerate in
the sense of Morse-Bott and Theorem~\ref{MB} implies that the degree
$\deg\pi$ is equal to $\pm\chi(M)$. Suppose that the homotopy class
$\ev_u^\natural[\phi]$ is not trivial. Then the energy
identity~\eqref{EnergyId} and the hypothesis~\eqref{NonPosive} imply
that this homotopy class does not contain $J$-tori. Hence, the space
of solutions $\pi^{-1}(0)$ is empty and the zero section of
$\Omega^{0,1}$ is regular. This implies that the degree $\deg\pi$
vanishes. Thus, we obtain that the Euler-Poincar\'e number $\chi(M)$
is equal to zero and the theorem is demonstrated.
\qed

\medskip
\noindent
{\bf The proof of Theorem~\ref{GS}$^{\mathbf{\prime}}$} follows along
similar lines, since the compactness theorem and Theorems~\ref{Folk}
and~\ref{MB} have analogues concerning the moduli space of
solutions within homology classes of mappings.

\medskip
We end the paper with explaining the proof of Lemma~\ref{moser}. 

\medskip
\noindent
{\bf Proof of Lemma~\ref{moser}.}
First, without loss of generality we can consider only maps that take
values in the connected component of the identity $\symp_0(M)$. 
Second, a given continuous map $\phi:\Sigma\to\symp_0(M)$ can be
$C^1$-approximated by a map $\hat\phi$ with values in $\diff_0(M)$
such that $\ev_u\circ\hat\phi$ is smooth for any $u\in M$. Indeed, one
can regard any map of a surface $\Sigma$ into the diffeomorphism group
as a map from the product $\Sigma\times M$ to $M$ and approximate it
by a smooth map with respect to the first variable. Finally, to obtain
an approximation with values in the symplectomorphism group we apply
to $\hat\phi$ the canonical retraction $\mathcal R$ of a
$C^1$-neighbourhood of $\symp_0(M)$ in the diffeomorphism group to
$\symp_{0}(M)$. We describe an explicit construction for the latter in
terms of the so-called Moser isotopy now. 

Let $\hat\phi$ be a diffeomorphism from $\diff_0(M)$ which is
$C^1$-close to $\symp(M)$ such that all forms
$$
\omega_t=\omega+t(\hat\phi^*\omega-\omega),\qquad t\in [0,1],
$$
are non-degenerate. Since $\hat\phi$ is homotopic to the identity, the
forms $\omega_t$'s are cohomologous. By Moser's isotopy
theorem~\cite[Section~3.2]{McDS0} there exists a canonical family of
diffeomorphisms such that $\psi_t^*\omega_t=\omega_0$ and
$\psi_0=\Id$. Clearly, the diffeomorphism $\psi_1\circ\hat\phi$
preserves the symplectic form $\omega$ and we define the retraction
$\mathcal R$ by the rule $\hat\phi\mapsto\psi_1\circ\hat\phi$. 

Now suppose that a diffeomorphism $\hat\phi$ depends smoothly on a
parameter $z\in\Sigma$ (in the sense that the map $\Sigma\times M\to
M$ is smooth), then so do the forms $\hat\phi^*\omega_t$'s. By the
construction of the Moser isotopy, the diffeomorphisms $\psi_t$ are
defined as solutions of certain differential equations and also depend
smoothly on the parameter $z$. Hence, so does the diffeomorphism
$\psi_1$ as well as the composition $\psi_1\circ\hat\phi$.
\qed

\end{document}